\documentstyle[12pt, oldlfont]{article}

\input{amssym.def}
\input{amssym.tex}

\begin{document}

\def\Ubf#1{{\baselineskip=0pt\vtop{\hbox{$#1$}\hbox{$\sim$}}}{}}
\def\ubf#1{{\baselineskip=0pt\vtop{\hbox{$#1$}\hbox{$\scriptscriptstyle\sim$}}}{}}
\def\R{{\Bbb R}}
\def\V{{\Bbb V}}
\def\N{{\Bbb N}}
\def\Q{{\Bbb Q}}
\def\O{{\cal O}}
\def\e{\epsilon}
\def\h{{\hat{X}}}
\def\B{{\cal B}}

\title{A universal Polish $G$-space}         
\author{Greg Hjorth \footnote{Research partially supported by NSF grant DMS 96-22977}}      
\date{\today}          
\maketitle

\leftskip 1in 

\noindent {\bf Abstract:} If $G$ is a Polish group, then the category of Polish 
$G$-spaces and continuous $G$-mappings has a universal object. 

\leftskip 0in 
\medskip




{\bf $\S$0. Preface} 

It is known from \cite{devries} that\\

0.1. Theorem (de Vries). Let $G$ be a separable locally compact metric group.  
Then the category  of continuous $G$-maps and separable complete metric $G$-spaces  
has a universal object.\\

In other words, there is a separable complete metric space $X$ on which $G$ acts 
continuously, and such that for any other separable completely metrizable 
$Y$ on which there is a continuous $G$ action we can find a continuous $\pi:X\rightarrow Y$ 
such that for all $x_0,x_1\in X$, $g\in G$
\[g\cdot x_0=x_1\Leftrightarrow g\cdot\pi(x_0)=\pi(x_1).\]

On the other hand in \cite{beckerkechris}:\\

0.2. Theorem (Becker, Kechris). Let $G$ be a Polish group -- that is to say a separable 
topological group admitting a complete metric. Then the category of Borel 
$G$-maps and separable complete metric $G$-spaces has a universal object.\\

While here we unify these with:\\

0.3. Theorem. Let $G$ be a Polish group. Then the category of continuous 
$G$-maps and separable complete metric $G$-spaces has a universal object.\\

The proof gives incidental information beyond the original goal. We may relax the 
condition that $G$ be Polish to simply being a separable metric group, and even in this 
case one has a Polish $G$-space that is universal for {\it all separable metric} $G$-spaces:\\

0.4. Theorem. Let $G$ be a separable metric group. Then there is Polish $X$ such that:

(i) $G$ acts continuously on $X$; 

\noindent and whenever $G$ acts continuously on some separable metric space $Y$ we may find 
$\pi:Y\rightarrow X$ so that 

(ii) $\pi$ is a continuous and one-to-one; 

(iii) $\pi$ is a $G$-mapping, in the sense that for all $g\in G, y\in Y (g\cdot\pi(y)=\pi(g\cdot y)).$\\

\newpage 

{\bf $\S$1 Definitions}\\

1.1. Definition. A topological space is said to be {\it Polish} if it separable 
and it has a complete metric compatible with its topology. A {\it Polish group} 
is a topological group that is Polish as a space.\\

Examples of Polish spaces are ${\Bbb R}$, ${\Bbb C}$, $\N^{\N}$, any compact metric 
space, and any separable manifold. More generally we might be forgiven for thinking that 
almost all topological spaces are Polish. The examples of Polish groups are equally 
broad, counting among their numbers $\l^2$, $c_0$, and more generally all separable Banach 
spaces, all separable Lie groups, and the symmetric group on the integers in the 
topology of pointwise convergence. When it becomes necessary to be specific about the 
topology I will write $(X, \tau)$ to denote a Polish space whose underlying set is 
$X$ and whose topology is $\tau$; the rest of the time we may as well follow the usual kind 
of convention and identify the Polish 
space with its underlying set.\\

1.2. Definition. If $G$ is a Polish group acting continuously on a Polish space 
$X$ then we say the $X$ is a {\it Polish $G$-space}. $X$ is said to be a {\it continuously universal 
Polish $G$-space} if for any other Polish $G$-space $Y$ there is a continuous $G$-injection  
$\pi:Y\rightarrow X$ -- so it is one-to-one, continuous, and for all 
$g\in G$, $y\in Y$ $\pi(g\cdot y)=g\cdot \pi(y).$ 
$X$ is a {\it Borel universal 
Polish $G$-space} if for any other Polish $G$-space $Y$ we may find Borel $G$-embedding 
from $X$ to $Y$.\\

1.3. Theorem (Birkhoff-Kakutani). Any Polish group has a compatible right invariant metric 
-- that is to say for all $g_0, g_1, h\in G$ $d(g_0,g_1)=d(g_0h, g_1h)$. 
(See \cite{hewittross}, 8.3)\\

This only promises that there will be a right invariant metric compatible with the 
topology on $G$; in fact Symm$(\N)$, the infinite permutation group, is an example 
of a Polish group that does not allow a {\it complete} right invariant metric. Letting 
$d_0(g,h)=d(g,h)(d(g,h)+1)^{-1}$ we obtain from 1.3 a right invariant metric bounded by 1.\\

1.4. Notation. Let $G$ be a Polish group and $d$ a right invariant metric bounded by 1. 
Then ${\cal L}(G,d)$ is the space of all $f:G\rightarrow [0,1]$ such that 
\[\forall g_1, g_2\in G(|f(g_1)-f(g_2)|\leq d(g_1, g_2).\]
We let $G$ act on this space by $(g\cdot f)(g_0)=f(g_0g)$ for all $g, g_0\in G$ and 
$f\in {\cal L}(G,d)$. 

For $A$ a set we can let ${\cal L}(G,d)^A$ be the space of 
all $\vec f:A\rightarrow {\cal L}(G,d)$. 
Here I will use $f_a$ as shorthand for the cumbersome $\vec f (a)$. 
We obtain an induced action on ${\cal L}(G,d)^A$ by taking the product of the 
actions on ${\cal L}(G,d)$, so that $(g\cdot \vec f)(a)=g\cdot(f_a)$.\\

Note that this is an action, since for all $\vec f\in{\cal L}(G,d)^A$, 
$a\in A$, $g, h, g_0\in G$ we have 
$(gh)\cdot f_a(g_0)=f_a(g_0gh)=h\cdot f_a(g_0g)=g\cdot(h\cdot f_a)(g_0)$. 

${\cal L}(G,d)$ has a natural topology obtained by viewing it as a closed subspace of 
$[0,1]^G$ in the product topology, and this in turn gives us a product topology on 
${\cal L}(G,d)^A$. Both these spaces are compact by Tychonov's theorem, 
and in the case that $A$ is countable they are separable, hence Polish. 
It is easily seen that the actions described above are then continuous, 
and so ${\cal L}(G,d)^A$ is a Polish $G$-space. 

While this space is a {\it Borel} universal $G$-space for $A$ countably infinite, it is 
known from \cite{meg} that in general we may find Polish $G$-spaces that allow 
no continuous embedding into any compact Polish $G$-space. The ultimate 
construction will consist in a carefully chosen subspace in a carefully chosen 
topology. Here we will in the end need:\\ 

1.5. Theorem (classical). Let $X$ be a Polish space and $X_0\subset X$ a $G_{\delta}$ 
subset. Then $X_0$ is Polish in the relative topology. (See \cite{kechris2}.)\\

\newpage 

{\bf $\S2$ The proof}

For the remainder of this section let $G$ be a Polish group, let $d$ be a right invariant 
metric compatible with its topology, and let $G_0$ be a countable dense subgroup of $G$. 
In what follows ${\Bbb Q}^{<\N}$ will denote the set of all finite sequences of rationals, and 
if $s\in {\Bbb Q}^{<\N}$ and $q\in{\Bbb Q}$, then $sq$ will denote the sequence of length one greater 
than $s$, extending $s$, and finishing with $q$ as its final value -- so that if 
$s=\langle q_0,q_1,...,q_n\rangle$ then $sq=\langle q_0,q_1,...,q_n,q\rangle$.\\

2.1. Definition. Let $X$ be the set of $\vec f\in {\cal L}(G,d)^{{\Bbb Q}^{<\N}}$ such that: 

\leftskip 1in

\noindent (1) for all $s\in \Q^{<\N}, g_0\in G_0, n, m, \in \N, 
q_0, q_1, q_2\in\Q\cap (0,1)$, $\epsilon\in \Q\cap (0,1)$, if $0<q_i\pm \epsilon<1$ 
for $i=0,1,2$, and if $t=sq_0q_1q_20mn$ then 
\[f_{t}(g_0)<q_1-\epsilon\vee f_s(g_0)\leq q_0+\epsilon;\]

\noindent (2) for all $s\in \Q^{<\N}, g_0\in G_0, n, m\in \N, 
q_0, q_1, q_2\in\Q\cap (0,1)$, $\epsilon\in \Q\cap (0,1)$, if $0<q_i\pm \epsilon<1$ 
for $i=0,1,2$, and if $u=sq_0q_1q_21mn$ then 
\[f_{u}(g_0)\geq q_2+\epsilon\vee f_{t}(g_0)\geq q_1-\epsilon;\] 

\noindent (3) for all $s\in \Q^{<\N}, g_0\in G_0, 
q_0\in\Q\cap (0,1)$, $\epsilon\in \Q\cap (0,1)$, if $f_s(g_0)<q_0$ then there is 
$q_2,q_1\in\Q$, $g_1\in G_0$ $n, m\in \N$ with 
$0<q_2<q_1<q_0<1$, $d(g_0, g_1)<\epsilon$, and for $u=sq_0q_1q_21mn$ 
\[f_{u}(g_1)<q_2.\]

\leftskip 0in

$X$ is a subset of ${\cal L}(G,d)^{{\Bbb Q}^{<\N}}$ which we will shortly equip with 
a rather different topology. On the other hand the action on $X$ will be as 
described at 1.4. In order for this to make sense:\\

2.2. Lemma. $X$ is $G$-invariant. 

Proof. We show this piecemeal, starting at (1), then (2), then (3). For these first two it suffices to 
see that (1) and (2) hold for {\it all } $g_0\in G$, not just $g_0\in G_0$. 

(1): Fix $s\in \Q^{<\N}, g_0\in G, n, m\in \N, 
q_0, q_1, q_2\in\Q\cap (0,1)$, $\epsilon\in \Q\cap (0,1)$, with  $0<q_i\pm \epsilon<1$ 
for $i=0,1,2$, and for $t=sq_0q_1q_20mn$
\[f_{t}(g_0)\geq q_1-\epsilon.\] 
Choose $\delta\in\Q^+$, $h_0\in G_0$ with $d(g_0, h_0)<\delta$, $0<q_i\pm(\epsilon+\delta)<1$ 
for $j=0,1,2$. Then we obtain that $f_{t}(h_0)\geq q_1-(\epsilon+\delta)$, 
so $f_s(h_0)\leq q_0+\epsilon+\delta$, by (1), and hence $f_s(g_0)\leq q_0+\epsilon+2\delta$. 
Letting $\delta$ tend to zero we have $f_s(g_0)\leq q_0+\epsilon$. 

(2): For $s\in \Q^{<\N}, g_0\in G, n, m\in \N, 
q_0, q_1, q_2\in\Q\cap (0,1)$, $\epsilon\in \Q\cap (0,1)$, 
$u=sq_0q_1q_21mn$ 
with  $0<q_i\pm \epsilon<1$ $i=0,1,2$, with 
\[f_{u}(g_0)< q_2+\epsilon.\] 
Now choosing $h_0$ in $G_0$ with $d(g_0, h_0)<\delta$, $0<q_i\pm(\epsilon+\delta)<1$ 
for $j=0,1,2$ we obtain essentially the same outcome as before in virtue of (2) from 2.1. 

(3): Fix $g_0\in G$, $\epsilon,  q_0\in (0,1)\cap \Q$, and suppose 
$f_s(g_0)<q_0$. Choose $\delta\in\Q\cap(0,\epsilon/3)$ so that $f_s(g_0)+\delta<q_0$. 
It suffices to show that there is an open neighbourhood of group elements 
$g_1$ that satisfy the conclusion of (3) for $g_0$ and some 
fixed $m,n, q_1, q_2$. 

Take $h_0\in G_0$ with $d(g_0,h_0)<\delta$. Then $f_s(h_0)<q_0$, so we may 
find $h_1\in G_0$ with $d(h_0, h_1)<\delta$ and $m, n, q_1, q_2$ so that 
for $u=sq_0q_1q_21mn$ we have 
\[f_{u}(h_1)<q_2.\] 
Then we may fix $\bar{\delta}<\delta$ so that 
\[f_{u}(h_1)<q_2-\bar{\delta}.\]
Then for all $g_1$ with $d(g_1, h_1)<\bar{\delta}$ we have $f_{u}(g_1)<q_2$ 
and $d(g_0, g_1)< d(g_0, h_0)+d(h_0, h_1)+d(h_1, g_1)<\delta+\delta+\delta<\epsilon$, 
as required. \hfill $\Box$\\

2.3 Definition of the topology on $X$: We take as a subasis for the topology $\tau$ all 
sets of the form $\{\vec f:f_s(g_0)<q_0\}$ for $g_0\in G_0$, $q_0\in\Q$,$s\in \Q^{<\N}$.\\

So at once $\tau$ is separable and $G$ acts continuously on $(X,\tau)$. Note that by the 
Megrelishvili result -- to the effect that there may be no {\it compact} universal Polish $G$-space -- 
$X$ cannot in general be a universal space in the product topology obtained from ${\cal L}(G,d)^{\Q^{<\N}}$.\\

2.4. Lemma. $(X,\tau)$ is Hausdorff. 

Proof. Fix $\vec f, \vec h\in X$ with $\vec f\neq \vec h$. 
Without loss of generality there is some $s\in \Q^{<\N}$, $q_0\in\Q$, $g_0\in G_0$, 
$\epsilon\in\Q^+$ with 
\[f_s(g_0)<q_0<q_0+2\epsilon<h_s(g_0),\]
\[0<q_0\pm\epsilon<1.\]
Then find $g_1\in G_0$, $q_2, q_1\in\Q$, $m,n$ as indicated by 2.1(3), 
so that for $t=sq_0q_1q_20mn$, $u=sq_0q_1q_21mn$, 
\[f_{u}(g_1)<q_2\]
\[\therefore f_{t}(g_1)\geq q_1,\]
by 2.1(2). On the other hand for $\delta<\epsilon$ sufficiently small $h_{t}(g_1)<q_1-\delta$ by (1). 
Then $\{\vec k:k_{t}(g_1)<q_2\}$ and $\{\vec k:k_{u}(g_1)<q_1\}$ 
are our disjoint open sets. \hfill $\Box$\\

2.5. Lemma. $(X,\tau)$ is regular. 

Proof. It suffices to show that if $\O$ is a subbasic open set and $\vec h\in\O$ 
then there is open $U\subset \O$ with $\vec h\in U$ and $\overline{U}$, the closure of 
$U$, included in $\O$. 

So suppose $\vec h\in \O=\{\vec f\in X: f_s(g_0)<\bar{q}_0\}$ for some 
$g_0\in G_0$, $\bar{q}_0\in\Q\cap(0,1)$. Then choose $q_0\in\Q\cap(0,\bar{q}_0)$ 
and $\epsilon\in\Q$ 
such that $h_s(g_0)<q_0<q_0+2\epsilon<\bar{q}_0$ and $q_0\pm\epsilon \in\Q\cap(0,1)$. 
By (3) we find $g_1\in G_0$ and $q_1, q_2, t, u$ as before with 
\[d(g_1, g_0)<\epsilon,\]
\[h_{u}(g_1)<q_2.\] 
Then for small $\delta<\epsilon$
\[\vec h\in\{\vec f: f_{u}(g_1)<q_2\}\subset \{\vec f: f_{t}(g_1)\geq q_1-\delta\},\]
which is closed and in turn included in 
\[\{\vec f: f_{s}(g_1)\leq q_0+\delta\}\subset \{\vec f: f_{s}(g_0)\leq q_0+2\epsilon\}
\subset\O,\]
as required. \hfill $\Box$\\

$(X,\tau)$ is a regular, separable, Hausdorff, and so is metrizable by Urysohn.\\

2.6. Definition. Let $d_0$ be some metric for $(X,\tau)$ and let 
$\h$ be the Cauchy completion of $(X,d_0)$. For each $s\in\Q^{<\N}$, $g_0\in G_0$, 
$q_0\in \Q$, let $\O(s, g_0, q_0)$ be an open set in $\h$ such that 
$\O(s, g_0, q_0)\cap X=\{\vec f : f_s(g_0)<q_0\}\cap X$.\\

2.7. Definition. Let $X_0\subset\h$ be the set of all $x\in\h$ such that 

\leftskip 1in

\noindent (4) for all $s\in \Q^{<\N}, g_0\in G_0, n,m\in \N, 
q_0, q_1, q_2\in\Q\cap (0,1)$, $\epsilon\in \Q\cap (0,1)$, $\delta\in\Q^+$, if $0<q_i\pm \epsilon<1$ 
for $i=0,1,2$, and if $t=sq_0q_1q_20mn$ then 
\[x\in\O(t, g_0, q_1-\epsilon)\vee x\in\O(s, g_0, q_0+\epsilon+\delta);\] 

\noindent (5) for all $s\in \Q^{<\N}, g_0\in G_0, n, m\in \N, 
q_0, q_1, q_2\in\Q\cap (0,1)$, $\epsilon\in \Q\cap (0,1)$, if $0<q_i\pm \epsilon<1$ 
for $i=0,1,2$, and if $u=sq_0q_1q_21mn$ then 
\[\neg(x\in\O(u,g_0,q_2+\epsilon))\vee \neg(x\in\O(t,g_0, q_1-\epsilon));\] 

\noindent (6) for all $s\in \Q^{<\N}, g_0\in G_0, 
q_0\in\Q\cap (0,1)$, $\epsilon\in \Q\cap (0,1)$, if $x\in\O(s,g_0,q_0)$ then there is 
$q_2,q_1\in\Q$, $g_1\in G_0$ $n, m\in \N$ with 
$0<q_2<q_1<q_0<1$, $d(g_0, g_1)<\epsilon$, and for $u=sq_0q_1q_21mn$ 
\[x\in\O(u,g_1,q_2);\]

\noindent (7) for all $s\in \Q^{<\N}, g_0, g_1\in G_0, 
q_0 \in\Q\cap (0,1)$, $\epsilon\in \Q\cap (0,1)$, if $d(g_0,g_1)<\epsilon$ then 
\[x\in\O(s, g_0, q_0)\Rightarrow x\in\O(s, g_1, q_0+\epsilon);\] 

\leftskip 0in 

\noindent which in particular implies 

\leftskip 1in 

\noindent (8) for all $s\in \Q^{<\N}, g_0\in G_0, 
q_0<q_1 \in\Q\cap (0,1)$, 
\[x\in\O(s, g_0, q_0)\Rightarrow x\in\O(s, g_0, q_1).\]\\

\leftskip 0in

2.8. Lemma. $X_0$ is $G_{\delta}$ in $\h$. 

Proof. In any metric space the closed sets are $G_{\delta}$, and thus all finite 
Boolean combinations of open sets are $G_{\delta}$; each of (4)-(7) corresponds to a countable 
intersection of $G_{\delta}$ sets, and is therefore $G_{\delta}$. \hfill $\Box$\\

2.9. Notation. For $x\in X_0$ let $\vec k^x\in{\cal L}(G_0,d)$ be given by $k^x_s(g_0)=$inf
$\{q\in\Q\cap(0,1):x\in\O(s,g_0,q_0)\}$ for $g_0\in G_0$ and $s\in\Q^{<\N}$. 

Since $G_0$ is dense, this extends canonically to a function $\vec f^x\in{\cal L}(G,d)$ by 
(7) above.\\

2.10. Lemma. $x\in X_0\Rightarrow \vec f^x\in X$. 

Proof. As before, (7) gives that by $\vec f^x\in{\cal L}(G,d)$ and $x\mapsto \vec f^x$ 
is welldefined. From (8) we have that for all $s\in\Q^{<\N}$, $g_0\in G_0$ 
\[x\in \O(s,g_0, q_0)\Leftrightarrow f_s^x(g_0)<q_0.\]
And then 2.1(1) follows from 2.7(4) and letting $\delta\rightarrow 0$, 
(2) from (5), and (3) from (6). \hfill $\Box$\\

On the other hand it is immediate that\\

2.11. Lemma. $X\subset X_0$.\\

And so by 2.11, 2.10, 1.5 and the definition of $\h$\\

2.12. Lemma. $X$ is a Polish $G$-space.\\

Our work will be over once we now verify:\\ 

2.13. Lemma. Let $Y$ be a Polish $G$-space. Then there is a continuous $G$-embedding 
$\pi:Y\rightarrow X$; moreover, we may choose $\pi$ to be open to its image, in the 
sense that if $U\subset Y$ is open the there is open $O\subset X$ such that 
$U=\{y\in Y:\pi(y)\in O\}$. 

Proof. Let $\B=(O_n)_{n\in \N}$ be a countable basis. For $y\in Y$ we now define 
$f^y_s$ uniformly in $y$ by induction on the length  of $s$. 

Case (i): For $s=\langle n \rangle$ some $n\in N$, we let $f_s^y(g_0)=$inf$\{d(1_G, g):gg_0\cdot y\in O_n\}$ 
for all $g_0\in G$, $y\in Y$ for which this infinum is defined; if $G\cdot y\cap O_n=\emptyset$ 
then let $f^y_s(g)=1$.

Case (ii): Suppose $s\in\Q^{<\N}$, and suppose that $O\subset Y$ is open and such that  
$f^y_s(g_0)=$inf$\{d(1_G, g):gg_0\cdot y\in O\}$ for all $y\in Y$ and $g_0\in G$ for which the 
infinum is defined (and $=1$ when undefined), 
and suppose that $t=sq_0q_1q_20mn$, $u=sq_0q_1q_21mn$ for some 
$q_2<q_1<q_0\in Q$, $m,n\in\N$. 

Case (iia): Suppose that 

\leftskip 1in 

\noindent (9) $\overline{O}_m\subset\bigcup\{g\cdot O:d(1_G,g^{-1})<q_0\}$, and 
that for all $z\in Y$, $g,h\in G$, 
\[d(1_G, g)<q_1+q_2\wedge g\cdot z\in O_n\wedge d(1_G,h)<q_1\Rightarrow h\cdot z \in O_m.\]

\leftskip 0in

Then let: 

\leftskip 0.8in 

$f^y_{t}(g_0)=$inf$\{d(1_G,g):gg_0\cdot y\in Y\setminus \overline{O}_m\}$ when defined 
($=1$ when the infinum is 
not defined),

$f^y_{u}(g_0)=$inf$\{d(1_G,g):gg_0\cdot y\in {O}_n\}$ when defined ($=1$ when not defined).  

\leftskip 0in 

Case (iib): Otherwise just let $f^y_{t}(g_0)=0$ and $f^y_{u}(g_0)=1$, all $y, g_0$. 

Case (iii): For $s$ not covered by the above cases, just let $f_s^y(g_0)=0$ all $y, g_0$. 

Claim I. For all $y\in Y$, $\vec f^y\in{\cal L}(G,d)$. 

Proof of claim: This follows by the fact that every $f^y_s$ is either a constant function, 
in which case it is trivial, or has the form $f^y_s(g_0)=$inf$\{d(1_G, g):gg_0\cdot y\in O\}$, 
when for $g_0, g_1$ with $d(g_0,g_1)<\delta$, and $d(1_G, g)<q_0$ with $gg_0\cdot y\in O$ 
witnessing $f^y_s(g_0)<q_0$, we have 
\[d(1_G, g_0g_1^{-1})<\delta,\]
\[d(1_G, g)=d(g_0g_1^{-1}, gg_0g_1^{-1})<q_0,\]
by right invariance, so that 
\[d(1_G, gg_0g_1^{-1})<\delta+q_0,\]
\[gg_0g_1^{-1}g_1\cdot y=gg_0\cdot y\in O\]
witnesses $f^y_s(g_1)<q_0+\delta$, as required. \hfill ($\Box$Claim I)

Claim II. For all $y\in Y$, $\vec f^y$ satisfies 2.1(1). 

Proof of claim: We check this for each $s\in \Q^{<\N}$; it is trivial 
except when we are in case (iia) above, so suppose that 
$f^y_s(g_0)=$inf$\{d(1_G, g):gg_0\cdot y\in O\}$ for all $g_0\in G$, 
$t=sq_0q_1q_20mn$, some 
$q_2<q_1<q_0\in Q$, $m,n\in\N$, 
$\overline{O}_m\subset\bigcup\{g\cdot O:d(1_G,g^{-1})<q_0\}=_{df}W$, and 
$f^y_{t}(g_0)=$inf$\{d(1_G,g):gg_0\cdot y\in Y\setminus \overline{O}_m\}$. 
Then if $f^y_{t}(g_0)\geq q_1-\epsilon$ some $\epsilon\in\Q\cap(0,q_1)$ it 
follows that in particular $g_0\cdot y\in \overline{O}_m\subset W$, and hence 
$f^y_{s}(g_0)<q_0<q_0+\epsilon$. \hfill ($\Box$Claim II) 

Claim III. For all $y\in Y$, $\vec f^y$ satisfies 2.1(2). 

Proof of claim: Again suppose $s\in\Q^{<\N}$ $O\subset Y$,  
$f^y_s(g_0)=$inf$\{d(1_G, g):gg_0\cdot y\in O\}$ for all $g_0\in G$, 
$t=sq_0q_1q_20mn$, $u=sq_0q_1q_21mn$ for some 
$q_2<q_1<q_0\in Q$, $m,n\in\N$, $q_j\pm \epsilon \in (0,1)\cap \Q$ for $j=0,1,2$,  
and that $O$, $O_n$, and $O_m$ are as described at (9) from case (iia). 
If $f^y_{u}(g_0)<q_2+\epsilon<q_1+q_2$ then fix $g$ with $d(1_G, g)<q_1+q_2$, 
$gg_0\cdot O_n$,  
and then note by (9) that for all $h$ with $d(1_G, h)<q_1$ we must have $hg_0\cdot y\in O_m$, 
and so $f^y_{t}(g_0)\geq q_1>q_1-\epsilon$ by definition of $\vec f ^y$. 

In the other cases 2.1(2) is immediate. \hfill ($\Box$Claim III)

Claim IV. For all $y\in Y$, $\vec f^y$ satisfies 2.1(3). 

So suppose that $O\subset Y$ is an open set, $s\in\Q^{<\N}$, 
and that for all $y$ and $g_0$ $f^y_s(g_0)=$inf$\{d(1_G, g):gg_0\cdot y\in O\}$ and that 
$\epsilon\in\Q^+$, $g_0\in G$, $y\in Y$, $q_0\in\Q\cap(0,1)$ with $f^y_s(g_0)<q_0$ we need to find 
suitable $q_1,$  $q_2$, $t$, $u$, and $g_1\in G_0$ with $f^y_s(g_1)<q_2$ and 
$d(g_0,g_1)<\epsilon$. Let $W=\bigcup\{g\cdot O:d(1_G,g^{-1})<q_0\}$ and by the regularity of the 
space choose basic open $O_m\subset W$ such that $g_0\cdot y\in O_m$ and $\overline{O}_m\subset W$. 
By the continuity of the action we may find $W_0\subset O_m$, $U\subset G$ both open so that 

\leftskip 1in 

\noindent (10) $1_G\in U$, $g_0\cdot y\in W_0$, $U\cdot W_0\subset O_m$, 

\leftskip 0in 

\noindent and then further open $O_n\subset W_0$, $U_0\subset U$ so that

\leftskip 1in 

\noindent (11) $1_G\in U_0$, $g_0\cdot y\in O_n$, $U_0\cdot O_n\subset W_0$.  

\leftskip 0in 

\noindent Then choose $\delta\in\Q\cap(0,q_0)$ so that for all $g\in G$

\leftskip 1in 

\noindent (12) $d(1_G, g)<\delta\vee d(1_G, g^{-1})<\delta \Rightarrow g\in U_0$. 

\leftskip 0in

\noindent Then we choose $q_2, q_1\in \Q$ so that $0<q_2<q_1<\delta/2<q_0$. We will be happy once we 
see that $t=sq_0q_1q_20mn$ and $u=sq_0q_1q_21mn$ satisfy (9). 

But if $z\in Y$, $g\in G$, $d(1_G, g)<q_2+q_1<\delta$ and $g\cdot z\in O_n$, 
then we may take $z_0=g\cdot z$, and note that by (12) $g^{-1}\in U_0$, and then by (11) 
$z=g^{-1}\cdot z_0\in W_0$, 
and so for all $h\in G$ with $d(1_G, h)<q_1<\delta$ we must have $h\cdot z\in O_m$ 
by (12) and (10) -- all 
as required for (9). So $f^y_{t}(g_0)=$inf$\{d(1_G,g):gg_0\cdot y\in Y\setminus \overline{O}_m\}$, 
$f^y_{u}(g_0)=$inf$\{d(1_G,g):gg_0\cdot y\in {O}_n\}$ because we are in case(iia). 
Since $g_0\cdot y\in O_n$ we get $f^y_{u}(g_0)=0<q_2$, which means we may fulfill 
(3) by taking $g_0=g_1$. \hfill ($\Box$Claim IV)

Claim V. $y\mapsto \vec f^y$ is a $G$-map. 

Proof of claim. Fix $s\in\Q^{<\N}$ so that $f_s^y(g_0)=$inf$\{d(1_G, g):gg_0\cdot y\in O\}$ all $y, g_0$.  
Then fix $y\in Y$ and $h, g_0\in G$. 
\[f^{h\cdot y}_s(g_0)=\mbox{inf}\{d(1_G, g):gg_0\cdot (h\cdot y)\in O\}\]
\[=\mbox{inf}\{d(1_G, g):g(g_0h)\cdot y\in O\}=f^y_s(g_0h)=h\cdot f^y_s(g_0).\]
\hfill ($\Box$Claim V)

Since $(O_n)_{n\in \N}$ is a basis for $Y$, $y\mapsto \vec f^y$ is one to one. 
By continuity of the $G$-action, the map is continuous and onto its image. \hfill $\Box$\\

So in other words:\\

2.14. Theorem. If $G$ is a Polish group, then there is a continuously universal 
Polish $G$-space.\\

In the course of the proof we only used that $G$ is a separable metric group to construct $X$, 
and for a given $Y$ we only need that it is separable metric $G$-space for the existence of 
a continuous $G$-embedding $\pi:Y\rightarrow X$.\\

\newpage

{\bf $\S$ 3 Concrete spaces}

The proof of 2.14 is abstract, showing us little more than the conclusion, 
and does not seem to give a $G$-space which is easy to understand. It is worth noting that 
that in specific cases one can obtain far more concrete universal spaces. 

The de Vries result gives that for $G$ locally compact we may take $L^2(G\times G)$ as 
our universal space. In general the Becker-Kechris construction gives a space only 
very marginally more satisfying than $X$ from $\S$2, but in the case of $G=$Symm$(\N)$ -- 
the group of all permutations of $\N$ with the topology of pointwise convergence -- 
one obtains a universal Borel $G$-space that is natural from the model theoretic perspective of 
\cite{hodges}. For each $n$ let $X_n=2^{\N^n}$ be the space of all functions from $\N^n$, the set of 
all length $n$ sequences of positive integers, to $\{0, 1\}$, with the discrete 
topology on $\{0,1\}$ and the induced product topology on on $X_n$; we then let Symm$(\N)$ act on 
$X_n$ by shift, so that for $f\in X_n$, $s=\langle k_1, k_2,...,k_n\rangle\in \N^n$
\[g\cdot f(s)=f(g^{-1}(k_1),...,g^{-1}(k_n)).\]
We then obtain a universal Symm$(\N)$-space by taking $\Pi_{n\in\N}X_n$. 

While this cannot be universal in the sense of continuous maps since it is a totally disconnected 
space, a variation of this and $\S$2 gives a somewhat improved continuously universal 
Polish Symm$(\N)$-space. For each $s\in\N^{<\N}$ with final value a strictly positive integer, $s(lh(s))=k$, 
let $X_s$ be the space of all functions from $\N^k$ to $\{0,1\}$, where $\{0,1\}$ is equipped with the 
topology under which $\{0\}$ is open but not closed -- so we take $X_s$ to be $X_k$ with the product 
of the Sierpinski topology on $\{0,1\}$. 
We then take the above shift action on $\Pi_{s\in\N^{<\N}}X_s$. 
Following the plan of $\S$2 we take $X\subset\Pi_{s\in\N^{<\N}}X_s$ to consist of the 
set of functions $\vec f$ such that for all $g\in$Symm$(\N)$:

\leftskip 1in 

\noindent (1) for $t=s01l$ with final value some stricly positive integer, 
$g\cdot f_t(1,2,..l)=1\Rightarrow g\cdot f_s(1,...,s(lh(s)))=0$; 

\noindent (2) for $u=s02m$, $t=s01l$, $g\cdot f_u(1,2,..m)=0\Rightarrow g\cdot f_t(1,...,l)=1$; 

\noindent (3) for $s$ with final value some stricly positive integer, if 
$g\cdot f_s(1,...,s(lh(s)))=0$ then there is some $u$ and $t$ as in (1) and (2) with 
$g\cdot f_u(1,...,m)=0$; 

\noindent and for $s$ not with final value a strictly positive integer, $f_s(1,...)$ is 
just constanly $0$ (or undefined). 

\leftskip 0in

This is a Polish Symm$(\N)$-space for the same reasons as before. It is a little more work 
to check that it is continuously universal; this requires the observation that 
${\cal B}_0=\{\{g\in$Symm$(\N):\forall i\leq k(g(i)=i)\}:k\in\N\}$ as neighbourhood basis, 
and that if ${\cal B}$ is a basis for the Polish Symm$(\N)$-space $Y$, 
then $\{U\cdot O:U\in{\cal B}_0, O\in{\cal B}\}$ is again a basis for $Y$. 

Given the failure of $\S$2 to provide a {\it natural} continuously universal Polish $G$-space 
we might ask for a better proof. Here I would imagine that any proof sufficing for all 
Polish groups is necessarily somewhat abstract and uninformative. Instead it might be 
more promising to look for natural continuously universal spaces only with specific 
Polish groups -- such as Symm$(\N)$, or $U_{\infty}$, the unitary group of $l^2$, or the 
group of all measure preserving transformations of the unit interval.

6363 MSB

Mathematics

UCLA

CA90095-1555

greg@math.ucla.edu

\end{document}